\documentclass[12pt]{amsart}

\usepackage{epsf}
\usepackage{amsmath,amssymb}
\usepackage{colordvi}

\usepackage[all]{xy}
\CompileMatrices
\newcommand{\xyinc}{\ar@{^{(}->}}
\numberwithin{equation}{section}
\theoremstyle{plain}
\newtheorem{thm}{Theorem}[section]
\newtheorem{prop}[thm]{Proposition}
\newtheorem{coro}[thm]{Corollary}

\theoremstyle{definition}
\newtheorem{defi}[thm]{Definition}
\newtheorem{rem}[thm]{Remark} %use \begin{rem}

\newcommand{\ten}{\mbox{\hspace*{-.5pt}\raisebox{1pt}{${\scriptstyle \otimes}$}
\hspace*{-4pt}}}
\newcommand{\st}{\mathrm{st}}          % for standarization
\newcommand{\Des}{\mathrm{Des}}
\newcommand{\GDes}{\mathrm{GDes}}             %for (proper) global descents
\newcommand{\AGDes}{\overline{\mathrm{GDes}}} %for all global descents

\newcommand{\setS}{\mathsf{S}}
\newcommand{\setT}{\mathsf{T}}
\newcommand{\setR}{\mathsf{R}}

\newcommand{\Q}{\mathbb{Q}}
\newcommand{\Z}{\mathbb{Z}}

\newcommand{\frakS}{\mathfrak{S}}

\newcommand{\calD}{\mathcal{D}}
\newcommand{\calF}{\mathcal{F}}

\newcommand{\calM}{\mathcal{M}}
\newcommand{\calQ}{\mathcal{Q}}
\newcommand{\calZ}{\mathcal{Z}}

\newcommand{\QSym}{\mathcal{Q}\mathit{Sym}}
\newcommand{\SSym}{\mathfrak{S}\mathit{Sym}}

\newcommand{\Sh}[1]{\mathfrak{S}^{#1}}

\newcounter{FNC}[page]
\def\newfootnote#1{{\addtocounter{FNC}{2}$^\fnsymbol{FNC}$%
     \let\thefootnote\relax\footnotetext{$^\fnsymbol{FNC}$#1}}}

\title[Hopf algebra of permutations]{Structure of the Malvenuto-Reutenauer\\
       Hopf algebra of permutations\\(Extended Abstract)} 
\author{Marcelo Aguiar}
\address{Department of Mathematics\\
         Texas A\&M University\\
         College Station\\
         Texas \ 77843\\
         USA}
\email{maguiar@math.tamu.edu}
\urladdr{http://www.math.tamu.edu/\~{}maguiar}

\author{Frank Sottile}
\address{Department of Mathematics\\
        University of Massachusetts\\
        Lederle Graduate Research Tower\\
        Amherst, Massachusetts, 01003\\
        USA}
\email{sottile@math.umass.edu}
\urladdr{http://www.math.umass.edu/\~{}sottile}

\thanks{Research of second author supported in part by NSF grant DMS-0070494}
\keywords{Hopf algebras, symmetric group, weak order, quasi-symmetric
        functions} 
\subjclass{Primary  05E05, 06A11, 16W30}
\date{11 March 2002}
\thanks{Extended abstract for FPSAC`02 in Melbourne.  Abridged version of~\cite{AS02}}
\begin{document}

\maketitle

\begin{center}
 \begin{minipage}[c]{380pt}
 \small
 {\sc Abstract.}
 We analyze the structure of the Malvenuto-Reutenauer
 Hopf algebra of permutations in detail.
 We give explicit formulas for its antipode, prove that it is a cofree
 coalgebra,  determine its  primitive elements and its coradical filtration and
 show that it decomposes  as a crossed product over the Hopf algebra of
 quasi-symmetric functions. 
 We also describe the structure constants of the multiplication as a
 certain number of facets of the permutahedron. 
 Our results reveal a close relationship between the structure of 
 this Hopf algebra and the weak order on the symmetric
 groups. 
\end{minipage}\vspace{20pt}\\

 \begin{minipage}[c]{380pt}
 \small
 {\sc R\'esum\'e.}
 On analyse la structure de l'alg\`ebre de Hopf
 de Malvenuto et Reutenauer en d\'etail. On donne des formules
 explicites pour son antipode, on prouve que c'est
 une coalg\`ebre colibre, on d\'etermine ses \'el\'ements
 primitifs et sa filtration coradicale et on montre
 qu'elle se d\'ecompose comme un produit crois\'e
 sur l'alg\`ebre de Hopf de fonctions quasi-sym\'etriques. On d\'ecrit aussi
 les constantes de structure de la multiplication comme
 un certain nombre de facettes du permuto\`edre. Nos r\'esultats
 mettent en \'evidence une forte relation entre la structure
 de cette alg\`ebre de Hopf et l'ordre faible dans les groupes sym\'etriques.
\end{minipage}
\end{center}
%%%%%%%%%%%%%%%%%%%%%%%%%%%%%%%%%%%%%%%%%%%%%%%%%%%%%%%%%%%%%%%%%%%%%%%%%%%%%%%%%%
\section*{Introduction}

Malvenuto~\cite{Malv} introduced the Hopf
algebra $\SSym$ of permutations, which has a linear
basis $\{\calF_u\mid u\in \frakS_n, n\geq0\}$ indexed by 
permutations in all symmetric groups $\frakS_n$.
The Hopf algebra $\SSym$ is non-commutative, non-cocommutative, self-dual,
and graded.
Among its sub- and quotient- Hopf algebras 
are many algebras central to algebraic combinatorics.
These include the algebra of symmetric functions~\cite{Mac,St99}, 
Gessel's algebra $\QSym$ of quasi-symmetric functions~\cite{Ges}, the algebra
of non-commutative symmetric functions~\cite{GKal}, the Loday-Ronco
algebra of planar binary trees~\cite{LR98}, Stembridge's algebra 
of peaks~\cite{Stem97}, the Billera-Liu algebra of Eulerian
enumeration~\cite{BilLiu}, and others.
The structure of these combinatorial Hopf algebras with respect to certain
distinguished bases has been an important theme in algebraic combinatorics,
with applications to the combinatorial problems these algebras were created
to study.
We give a detailed understanding of the
structure of $\SSym$, both in algebraic and in combinatorial terms.

Our main tool is a new basis
$\{\calM_u\mid u\in \frakS_n, \ n\geq0\}$ for $\SSym$ related to its original
basis by M\"obius inversion on the weak order of the
symmetric groups. 
These bases $\{\calM_u\}$ and $\{\calF_u\}$ are analogous to the
monomial and fundamental basis of $\QSym$,
which are related via M\"obius inversion on their index sets, the Boolean
posets $\calQ_n$.

We give enumerative-combinatorial descriptions of the product, coproduct,
and antipode of $\SSym$ with respect to the basis
$\{\calM_u\}$. 
For example, the coproduct is obtained by splitting a
permutation at certain special positions that we call global descents. 
Descents and global descents
are left adjoint and right adjoint to a natural map $\calQ_n\to\frakS_n$.

The structure constants for the product with respect to the basis
$\{\calM_u\}$ are non-negative integers with the following
geometric-combinatorial description.
The $1$-skeleton of the permutahedron $\Pi_{n-1}$ is the Hasse diagram of the
weak order on $\frakS_n$.
The facets of the permutahedron are canonically isomorphic to products of lower
dimensional permutahedra. Say that
a facet isomorphic to $\Pi_{p-1}\times\Pi_{q-1}$ has type $(p,q)$.
Given $u\in\frakS_p$ and $v\in\frakS_q$, such a facet has a distinguished
vertex corresponding to $(u,v)$ under the canonical isomorphism.
Then, for $w\in\frakS_{p+q}$, the coefficient of
$\calM_w$ in the product $\calM_u\cdot\calM_v$ is the number of facets of the 
permutahedron $\Pi_{p+q-1}$ of type $(p,q)$ with the property that the
distinguished vertex is below $w$ and closer to $w$ than
to any other vertex in the facet.

We also give explicit formulas for the antipode with
respect to both bases. 
The structure constants with respect to the basis $\{\calM_u\}$
have constant sign, as in the case of $\QSym$.
The situation is more complicated for the basis $\{\calF_u\}$,
which may explain why no such explicit formulas were previously known.
\smallskip

Elucidating the elementary structure of $\SSym$ with respect to the
basis reveals further algebraic structures. 
For example, $\SSym$ is a cofree graded coalgebra. 
A consequence
is that the coradical filtration of $\SSym$ 
(which encapsulates the complexity of iterated coproducts) is the algebraic
counterpart of a filtration of the symmetric groups by certain lower order
ideals. In particular, the space of primitive elements is spanned
by the set $\{\calM_u\mid u \text{ has no global descents}\}$. 
Cofreenes was shown by Poirier and Reutenauer~\cite{PR95}, in dual form,
through the introduction of a different basis. 
The study of primitive elements was pursued from this point of view 
by Duchamp, Hivert, and Thibon~\cite{DHT01}. 

There is a well-known morphism of Hopf algebras $\SSym\to\QSym$ that maps
one fundamental basis onto the other, by associating to a permutation $u$ its
descent set $\Des(u)$. 
We describe this map in terms of the 
bases  $\{\calM_u\}$ and $\{M_\alpha\}$.

Lastly, $\SSym$ decomposes as a crossed
product over $\QSym$. This construction from the theory of Hopf algebras
is a generalization of the notion of group extensions. We provide a
combinatorial description for the Hopf kernel of the map $\SSym\to\QSym$.

These results are expanded on and proven in the manuscript~\cite{AS02} of the 
same name.
For a background on quasi-symmetric functions, see~\cite[\S 7.19]{St99}, for
Hopf Algebras, we recommend the book of Montgomery~\cite{Mo93a}.
We also recommend the papers~\cite{PR95} of Poirier and Reutenauer
and~\cite{DHT01} of Duchamp, Thibon, and Hivert, who studied this same
Hopf algebra of permutations from a different perspective, the latter under
the name `free quasi-symmetric functions'.

We thank Swapneel Mahajan, Nantel Bergeron, and the referees
for helpful comments.
%%%%%%%%%%%%%%%%%%%%%%%%%%%%%%%%%%%%%%%%%%%%%%%%%%%%%%%%%%%%%%%%%%%%%%%%%

\section{Essential definitions}

\subsection{Quasi-symmetric functions}

Gessel~\cite{Ges} introduced the algebra $\QSym$ of quasi-symmetric functions 
as the natural target for Stanley's $P$-partition generating function.
Subsequent work has shown its centrality, even universality, for generating
functions in algebraic combinatorics~\cite{Eh96,BMSW00,Ag01}.

A sequence $\alpha=(\alpha_1,\alpha_2,\ldots,\alpha_k)$ of positive integers
is a {\it composition of $n$} if $\sum_i\alpha_i=n$.
Compositions of $n$ correspond to subsets  of $[n{-}1]$ as follows
 \[
    \alpha\ =\ (\alpha_1,\alpha_2,\ldots,\alpha_k)
      \quad\longleftrightarrow\quad
    I(\alpha)\ =\ 
    \{\alpha_1,\alpha_1+\alpha_2,\dotsc,\alpha_1+\dotsb+\alpha_{k-1}\}\ .
 \]
These subsets of $[n{-}1]$ (and thus compositions of $n$) form the Boolean
poset $\calQ_n$, and the induced order relation on compositions is called 
{\it refinement}.

For $\setS\subseteq [n{-}1]$, the {\it fundamental quasi-symmetric} 
function $F_{\setS,n}$ is 
 \[
   F_{\setS,n}\ :=\ \sum_{\substack{j_1\leq\dotsb\leq j_n\\
                            i\in\setS\Rightarrow j_i<j_{i+1}}}
              x_{j_1}x_{j_2}\dotsc x_{j_n}\ .
 \]
These form a basis for %the graded $\Q$-vector space %of quasi-symmetric functions.
$\QSym$.
Another basis is provided by the monomial quasi-symmetric functions $M_\alpha$,
which are indexed by compositions $\alpha=(\alpha_1,\ldots,\alpha_k)$
 \[
   M_\alpha\ :=\ \sum_{i_1<i_2<\dotsb<i_k}
            x_{i_1}^{\alpha_1}x_{i_2}^{\alpha_2}\dotsb x_{i_k}^{\alpha_k}\ .
 \]
These two bases are related via M\"obius inversion on the Boolean poset
$\calQ_n$.
 \[
    F_\alpha\ =\ \sum_{\alpha\leq\beta} M_\beta
    \qquad\mbox{ and }\qquad
    M_\alpha\ =\ \sum_{\alpha\leq\beta}(-1)^{c(\beta)-c(\alpha)}F_\beta\,,
 \]
where $c((\alpha_1,\ldots,\alpha_k))=k$ and 
$(-1)^{c(\beta)-c(\alpha)}$ is the M\"obius function of $\calQ_n$.

The product of these $M_\alpha$ is given by quasi-shuffles of their 
indices~\cite[Lemma 3.3]{Eh96}. 
A \emph{quasi-shuffle} of compositions $\alpha$ and $\beta$ is a shuffle of
the components of $\alpha$ and $\beta$, where in addition we may replace
any number of pairs of consecutive components $(\alpha_i,\beta_j)$ in the
shuffle by $\alpha_i+\beta_j$.
Then we have
 \begin{equation}\label{E:prodqsym}
   M_\alpha \cdot M_\beta\ =\ \sum_\gamma M_\gamma\,,
 \end{equation}
where the sum is over all quasi-shuffles $\gamma$ of the compositions $\alpha$
and $\beta$. For instance,
 \begin{equation}\label{E:prod-ex}
 M_{(2)} \cdot M_{(1,1)}  \ =\ 
   M_{(1,1,2)}+M_{(1,2,1)}+M_{(2,1,1)}+M_{(1,3)}+M_{(3,1)}\,.
 \end{equation}
The unit element $1=M_{(\ )}$ is indexed by the empty composition.

Let $X$ and $Y$ be ordered alphabets with $X<Y$ their disjoint union
ordered as indicated.
Substitution $f(X)\mapsto f(X < Y)$ induces a coproduct
$\Delta\colon \QSym\to\QSym\ten\QSym$ whose action on a monomial function is
as follows.
\begin{equation}\label{E:copqsym}
  \Delta\bigl(M_{(\alpha_1,\ldots,\alpha_k)}\bigr)\ =\,\ \sum_{p=0}^k
  M_{(\alpha_1,\ldots,\alpha_p)}\ten M_{(\alpha_{p+1},\ldots,\alpha_k)}\,.
\end{equation}
For instance, 
$\Delta(M_{(2,1)})=1\ten M_{(2,1)}+M_{(2)}\ten M_{(1)}+M_{(2,1)}\ten 1$.

These definitions give $\QSym$ the structure of a graded, connected Hopf
algebra.
The degree $n$ component is spanned by those $M_\alpha$ where $\alpha$ is a
composition of $n$.
It is connected, as its degree $0$ component is 1-dimensional, 
and it is a Hopf algebra.
An explicit formula for the antipode was given by 
Malvenuto~\cite[corollaire 4.20]{Malv} and 
Ehrenborg~\cite[Proposition 3.4]{Eh96} 
 \begin{equation}\label{E:Q-antipode}
   S(M_\alpha)\ =\
   (-1)^{c(\alpha)}\sum_{\beta\leq\alpha}M_{\widetilde{\beta}}\,.
 \end{equation}
Here, if $\beta=(\beta_1,\beta_2,\ldots,\beta_t)$ then 
$\widetilde{\beta}$ is $\beta$ written in reverse order 
$(\beta_t,\ldots,\beta_2,\beta_1)$.

%%%%%%%%%%%%%%%%%%%%%%%%%%%%%%%%%%%%%%%%%%%%%%%%%%%%%%%%%%%%%%%%%%%%%%%%%%%
\subsection{The Hopf algebra of permutations}

Let $\SSym$ be the graded $\Q$-vector space with {\it fundamental basis}
$\{\calF_u\mid u\in \frakS_n, n\geq 0\}$, graded by $n$. 
$\SSym$ has a graded Hopf algebra structure first considered in
Malvenuto's thesis~\cite[\S 5.2]{Malv} and in her work with
Reutenauer~\cite{MR95}. 
Write $1$ for the basis element of degree $0$.

The product of two basis elements is obtained by shuffling the corresponding
permutations, as in the following example.
 \begin{eqnarray*}
  \calF_{\Blue{12}}\cdot\calF_{\Brown{312}} &=& \ \ \ 
    \calF_{\Blue{12}\Brown{534}}\,+\,
    \calF_{\Blue{1}\Brown{5}\Blue{2}\Brown{34}}\,
    +\,\calF_{\Blue{1}\Brown{53}\Blue{2}\Brown{4}}\,+\,
    \calF_{\Blue{1}\Brown{534}\Blue{2}}\,+\,
    \calF_{\Brown{5}\Blue{12}\Brown{34}}\\
    && +\,\calF_{\Brown{5}\Blue{1}\Brown{3}\Blue{2}\Brown{4}}\,+\,
    \calF_{\Brown{5}\Blue{1}\Brown{34}\Blue{2}}\,+\,
    \calF_{\Brown{53}\Blue{12}\Brown{4}}\,
    +\,\calF_{\Brown{53}\Blue{1}\Brown{4}\Blue{2}}\,+\,
    \calF_{\Brown{534}\Blue{12}}\,.
 \end{eqnarray*}
More precisely, for $p,q>0$, set
 \[
   \Sh{(p,q)} \ :=\
              \{\zeta\in \frakS_{p+q}\mid \zeta \mbox{ has at most one
              descent, at position $p$}\}\,.
 \]
This is the  collection of minimal (in length) representatives of left cosets
of the Young or parabolic subgroup $\frakS_p\times\frakS_q$ in $\frakS_{p+q}$,
called {\it Grassmannian permutations}. 
With these definitions, we describe the product.
For $u\in \frakS_p$ and $v\in \frakS_q$, set
 \begin{equation}\label{E:prod-fundamental}
  \calF_u \cdot \calF_v\ =\ \sum_{\zeta\in \Sh{(p,q)}}
                            \calF_{(u\times v)\cdot\zeta^{-1}}. 
\end{equation}
This endows $\SSym$ with the structure of a graded algebra with unit 1.

The algebra $\SSym$ is also a graded coalgebra with coproduct given by all
ways of splitting a permutation.
For a sequence $(a_1,\ldots,a_p)$ of distinct  integers, let
its {\it standard permutation}\newfootnote{Some authors call this 
flattening.} $\st(a_1,\ldots,a_p)\in\frakS_p$
be the permutation $u$ defined by  
\begin{equation}\label{E:st}
   u_i<u_j \iff a_i<a_j.
\end{equation}
For instance, $\st(625)=312$.
The coproduct $\Delta\colon\SSym\to\SSym\,\ten\SSym$ is defined by
\begin{equation}\label{E:cop-malvenuto}
  \Delta(\calF_u)\ =\ \sum_{p=0}^n \calF_{\st(u_1,\,\ldots,\,u_p)}\ten
                                   \calF_{\st(u_{p+1},\,\ldots,\,u_n)}\,,
\end{equation}
when $u\in\frakS_n$.
For instance, $\Delta(\calF_{42531})$ is 
\[
  1\ten\calF_{42531} +\calF_{1}\ten\calF_{2431} + \calF_{21}\ten\calF_{321} 
  +\calF_{213}\ten\calF_{21}+\calF_{3142}\ten\calF_{1} +\calF_{42531}\ten 1\,.
\]
$\SSym$ is a graded connected Hopf algebra~\cite[th\'eor\`eme 5.4]{Malv}.

This Hopf algebra $\SSym$ has been an object of
recent interest~\cite{Re93, MR95, PR95, RP01, DHT, DHT01, LR98, LR01}.
We remark that sometimes it is the dual Hopf algebra that is considered.
To compare results, one may use that $\SSym$ is
self-dual under the map $\calF_u\mapsto\calF_{u^{-1}}^*$, where
$\calF_{u^{-1}}^*$ is the element of the dual basis that is dual to 
$\calF_{u^{-1}}$.

To define the \emph{monomial} basis $\{\calM_u\}$ for $\SSym$
(in analogy to the basis $\{M_\alpha\}$ of $\QSym$),
we use the weak order on the symmetric groups $\frakS_n$.
Let $\ell(n)$ count the inversions $\{i<j \mid u_i>u_j\}$ of a permutation
$u$. 
The {\it weak order} on $\frakS_n$ is defined by 
 \[
   u\leq v\iff \exists w\in\frakS_n 
      \text{ such that }v\ =\ wu \text{ and } \ell(v)=\ell(w)+\ell(u)\,.
 \]
The cover relation $u\lessdot v$ occurs precisely when $v$ is obtained from 
$u$ by transposing a
pair of consecutive {\sl values} of $u$; a pair $(u_i,u_j)$ such that $i<j$ and
$u_j=u_i+1$. 
The maximum element of $\frakS_n$ is $\omega_n=(n,\ldots,2,1)$. 
\begin{figure}[htb]
  $$  \epsfxsize=2.8in\epsfbox{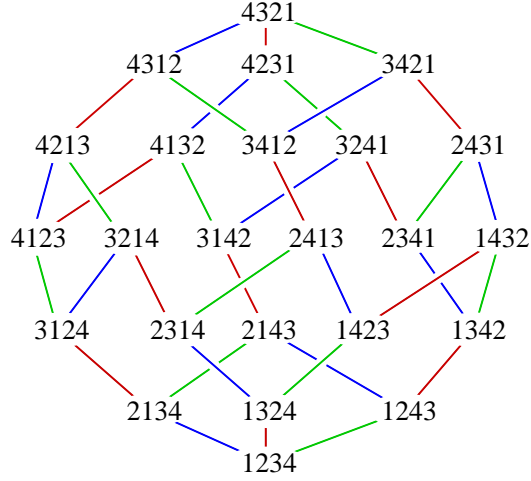}$$
  \caption{The weak order on $\frakS_4$\label{F:S4}}
\end{figure}
Figure~\ref{F:S4} shows the weak order on $\frakS_4$.

For each $n\geq 0$ and $u\in \frakS_n$, define
 \begin{equation}\label{E:def-monomial}
   \calM_u\ :=\ \sum_{u\leq v} \mu_{\frakS_n}(u,v)\cdot \calF_v\,,
 \end{equation}
where $\mu_{\frakS_n}(\cdot,\cdot)$ is the M\"obius function of 
the weak order in $\frakS_n$. 
By M\"obius inversion, 
 \begin{equation}\label{E:fun-mon}
   \calF_u\ :=\ \sum_{u\leq v} \calM_v\,,
 \end{equation}
so these elements $\calM_u$ indeed form a basis of $\SSym$. 
For instance,
 \[
   \calM_{4123}\ =\ \calF_{4123}-\calF_{4132}-\calF_{4213}+\calF_{4321}\,.
 \]
%

%%%%%%%%%%%%%%%%%%%%%%%%%%%%%%%%%%%%%%%%%%%%%%%%%%%%%%%%%%%%%%%%%%%%%%%%%%%%%%%
\subsection{The descent map $\calD\colon\SSym\to\QSym$}\label{S:monomialbasis} 

The descent set $\Des(u)$ of a permutation $u\in\frakS_n$ is the
subset of $[n{-}1]$ recording the descents of $u$
 \begin{equation}\label{E:defdescents}
   \Des(u)\ :=\ \{p\in[n{-}1]\mid u_p>u_{p+1}\}\,.
 \end{equation}
Thus $\Des(46\mspace{2mu}5\mspace{1mu}128\mspace{2mu}37)=\{2,3,6\}$.
Malvenuto~\cite[th\'eor\`emes 5.12, 5.13, and 5.18]{Malv} shows that 
there is a morphism of Hopf algebras
 \begin{equation} \label{E:descentmap}
   \begin{array}{rcrcl}
     \calD &:& \SSym&\longrightarrow& \QSym\\
           & &\calF_u&\longmapsto& F_{\Des(u)}\rule{0pt}{14pt}
   \end{array}
 \end{equation}
%
%(This is equivalent to Theorem 3.3 in~\cite{MR95}.)
%We show that $\calM_u$ maps either to $M_{\Des(u)}$ or to 0 under the map
%$\SSym\to\QSym$.
%%%%%%%%%%%%%%%%%%%%%%%%%%%%%%%%%%%%%%%%%%%%%%%%%%%%%%%%%%%%%%%%%%%%%%%%%%%

\subsection{Galois connections $\frakS_n\leftrightarrows\calQ_n$}\label{S:descents}
Underlying our results are combinatorial facts
concerning the lattices $\frakS_n$ and $\calQ_n$.
We describe two important conceptual facts.
For a subset $\setS\subseteq[n{-}1]$, let $\frakS_\setS\subseteq\frakS_n$ be the
parabolic subgroup
 \[
    \frakS_{\alpha_1}\times\frakS_{\alpha_2}\times\dotsb\times\frakS_{\alpha_k}\,,
 \]
where $\alpha=(\alpha_1,\dotsc,\alpha_k)$ is the composition of $n$ such that
$I(\alpha)=\setS$. 
For a subset $\setS\subseteq[n{-}1]$, let $Z(\setS)\in\frakS_n$ be the
maximal permutation with descent set $\setS$.

%\begin{defi}
 A \emph{Galois connection} between posets $P$ and $Q$ is a pair $(f,g)$ of
 order preserving maps $f\colon P\to Q$ and $g\colon Q\to P$ such that for any
 $x\in P$ and $y\in Q$,
 \begin{equation}\label{E:galois}
   f(x)\ \leq\ y \ \iff\  x\ \leq\ g(y)\,.
 \end{equation}
 Equivalently, the map $f$ is left adjoint to the map $g$.
%\end{defi}

\begin{prop}\label{P:galois}
 The pair of maps $(\Des,Z):\frakS_n\rightleftarrows\calQ_n$
 is a Galois connection.
\end{prop}

\begin{figure}[htb]
\[
  \begin{picture}(180,92)
   \put(0,0){
    \begin{picture}(70,92)
                        \put(25,83){321}
    \put(11,68){\line(1,1){10}}  \put(55,68){\line(-1,1){10}}
    \put( 0,55){312}             \put(50,55){231}
    \put( 8,40){\line(0,1){10}}  \put(60,40){\line(0,1){10}}
    \put( 0,27){213}             \put(50,27){132}
    \put(21,12){\line(-1,1){10}} \put(45,12){\line(1,1){10}}
                        \put(24,0){123}
    \end{picture}
   }
   \put(85,40){$\stackrel{\Des}{\longmapsto}$}

   \put(85,40){$\stackrel{\Des}{\longmapsto}$}
   \put(113,0){
    \begin{picture}(70,92)(1,0)
                        \put(19.5,83){\{1,2\}}
    \put(11,68){\line(1,1){10}}  \put(55,68){\line(-1,1){10}}
    \put(-1,55){\{1\}}             \put(50,55){\{2\}}
    \put( 7,40){\line(0,1){10}}  \put(58,40){\line(0,1){10}}
    \put( 9,40){\line(0,1){10}}  \put(60,40){\line(0,1){10}}
    \put(-1,27){\{1\}}             \put(50,27){\{2\}}
    \put(24,11){\line(-1,1){10}} \put(42,11){\line(1,1){10}}
                        \put(30,0){$\emptyset$}
    \end{picture}
   }
   \end{picture}
   \qquad\qquad
  \begin{picture}(170,92)
   \put(0,0){
    \begin{picture}(70,92)(1,0)
                        \put(19.5,69){\{1,2\}}
    \put(11,54){\line(1,1){10}}  \put(55,54){\line(-1,1){10}}
    \put(-1,41){\{1\}}             \put(50,41){\{2\}}
    \put(24,24){\line(-1,1){10}} \put(42,24){\line(1,1){10}}
                        \put(30,13){$\emptyset$}
    \end{picture}
   }
   \put(80,40){$\stackrel{Z}{\longmapsto}$}

   \put(103,0){
    \begin{picture}(70,92)
                        \put(25,69){321}
    \put(11,54){\line(1,1){10}}  \put(55,54){\line(-1,1){10}}
    \put( 2,41){312}             \put(53,41){231}
    \put(22,24){\line(-1,1){11}} \put(43,24){\line(1,1){11}}
                        \put(24,13){123}
    \end{picture}
   }
  \end{picture}
\]
 \caption{The Galois connection $\frakS_3\rightleftarrows\calQ_3$\label{F:galois}}
\end{figure}
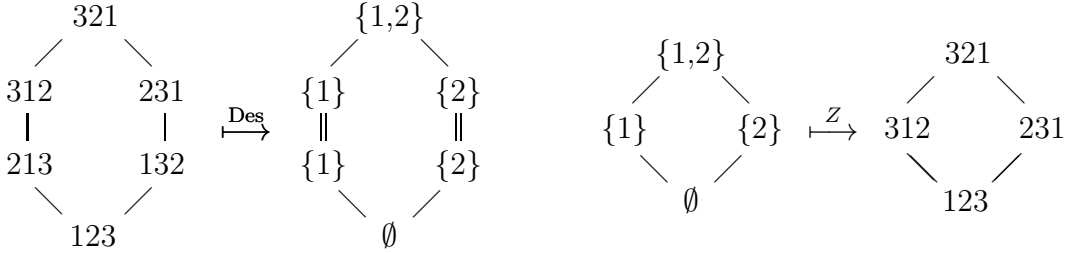

This Galois connection is why the monomial basis  of
$\SSym$ is analogous to that of $\QSym$, and is why we consider
the weak order on $\frakS_n$.
The connection between the monomial bases of these two algebras
will be elucidated in Theorem~\ref{T:map-monomial}.

 A permutation $u\in\frakS_n$ has a \emph{global descent} at a position
 $p\in[n{-}1]$ if 
  \[ 
    i\leq p<j \ \Longrightarrow\ u_i>u_j\,.
  \]
 Equivalently, if $\{u_1,\ldots,u_p\}=\{n,n{-}1,\ldots,n{-}p{+}1\}$.
 Let $\GDes(u)\subseteq[n{-}1]$ be the set of global descents of $u$.
 Note that $\GDes(u)\subseteq \Des(u)$, but these are not equal in
 general.

The notion of global descents is a very natural companion of that of
(ordinary) descents, in that the map
$\GDes\colon\frakS_n\to\calQ_n$
is {\em right} adjoint to $Z\colon\calQ_n\to\frakS_n$.
\begin{prop}\label{P:galoisglobal}
 The pair of maps $(Z,\GDes):\calQ_n\rightleftarrows\frakS_n$
 is a Galois connection.
\end{prop}

%%%%%%%%%%%%%%%%%%%%%%%%%%%%%%%%%%%%%%%%%%%%%%%%%%%%%%%%%%%%%%%%%%%%%%%%%%%%%%%%

\section{Elementary Algebraic Structure of $\SSym$}

\subsection{The coproduct of $\SSym$}\label{S:coproduct}

The coproduct of $\SSym$~\eqref{E:cop-malvenuto} takes a simple form on
the monomial basis.
For a permutation $u\in\frakS_n$, define $\AGDes(u)$ to be
$\GDes(u)\cup\{0,n\}$.

\begin{thm}\label{T:cop-monomial}
 Let $u\in \frakS_n$.
 Then
 \begin{equation}\label{E:cop-monomial}
    \Delta(\calM_u)\ =\ \sum_{p\in\AGDes(u)}
    \calM_{\st(u_1,\ldots,u_p)}\ten\calM_{\st(u_{p+1},\ldots,u_n)}\,.
 \end{equation}
\end{thm}

%%%%%%%%%%%%%%%%%%%%%%%%%%%%%%%%%%%%%%%%%%%%%%%%%%%%%%%%%%%%%%%%%%%%%%%%%
\subsection{The product of $\SSym$}\label{S:product}
The product of $\SSym$ in terms of its monomial basis 
has non-negative structure constants, which we describe.
For instance, 
 \begin{multline}\label{E:prod-M}
  \calM_{12} \cdot\calM_{21}\ =\ \calM_{4312}+\calM_{4231}
    +\calM_{3421}+\calM_{4123}+\calM_{2341}\\
    +\calM_{1243}+\calM_{1423}+\calM_{1342}+3\calM_{1432}+
     2\calM_{2431}+2\calM_{4132}\,. 
 \end{multline}
First, for a Grassmannian permutation $\zeta\in\frakS^{(p,q)}$
(a left coset representative of 
$\frakS_p\times \frakS_q$ in $\frakS_{p+q}$), consider the map corresponding
to the {\it right coset} of $\zeta^{-1}$.
 \[
   \rho_\zeta\ \colon\ \frakS_p\times\frakS_q\ \to\ \frakS_{p+q}\,,
    \qquad \rho_\zeta(u,v)\ :=\ (u\times v)\cdot\zeta^{-1}\,.
 \]
This order-preserving map is injective and its image is an interval in
$\frakS_{p+q}$.
For $u\in\frakS_p$, $v\in\frakS_q$ and $w\in\frakS_{p+q}$,
define $A^w_{u,v}\subseteq\Sh{(p,q)}$ to be 
 \begin{equation}\label{E:defalpha}
   A^w_{u,v}\ =\ \bigl\{\zeta \in \Sh{(p,q)}\mid
                   (u,v)=\max\rho_\zeta^{-1}[1,w]\bigr\}\,,
 \end{equation}
where $[w,w']:=\{w''\mid w\leq w''\leq w'\}$ denotes
the interval between $w$ and $w'$.
This set has another description as the set of those $\zeta\in\Sh{(p,q)}$ 
satisfying
 \begin{equation}\label{E:A-def}
    \begin{array}{ll}
       (i) & (u\times v)\cdot\zeta^{-1} \leq w, \text{ and}\\
      (ii) & \text{if $u\leq u'$ and $v\leq v'$ satisfy
             $(u'\times v')\cdot\zeta^{-1}\leq w$,}\\
           & \text{then $u=u'$ and $v=v'$.}
     \end{array}
 \end{equation}
Set $\alpha^w_{u,v}:=\#A^w_{u,v}$.

\begin{thm}\label{T:prod-monomial}
  For any $u\in\frakS_p$ and $v\in\frakS_q$, we have
 \begin{equation}\label{E:prod-monomial}
    \calM_u\cdot\calM_v\ =\ \sum_{w\in\frakS_{p+q}}\alpha^w_{u,v}\,\calM_w\,.
 \end{equation}
\end{thm}

For instance, in~\eqref{E:prod-M} the coefficient of $\calM_{2431}$ in 
$\calM_{12}\cdot\calM_{21}$ is 2 because among the six permutations
in $\Sh{(2,2)}$, %
 \[
    1234,\ 1324,\ 1423,\ 2314,\ 2413,\ 3412\,,
 \]
only the first two satisfy conditions $(i)$ and $(ii)$ of~\eqref{E:A-def}.

\begin{figure}[htb]
 \[
   \setlength{\unitlength}{0.7pt}%
   \begin{picture}(230,215)(-10,0)
    \put(0,0){\epsfysize=147.525pt\epsffile{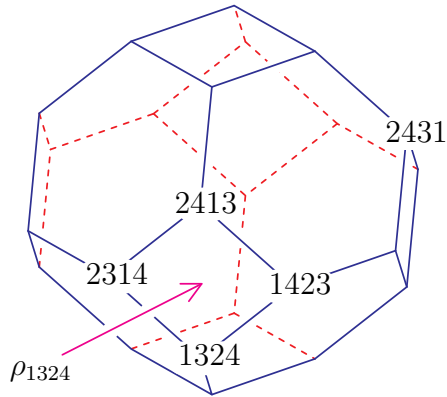}}
    \put(193,138){2431}
    \put(79,99.5){2413}
    \put(31,61){2314}  \put(130,56){1423}
    \put(81,17.5){1324}
    \put(-10,12){$\rho_{1324}$}
   \end{picture}
 \]
  \caption{The facet $\rho_{1324}$ of type $(2,2)$ and
    $w=2431$\label{F:alphaperm}.} \end{figure}

The structure constants $\alpha^w_{u,v}$ admit a geometric-combinatorial
description in terms of the permutahedron. 
The vertices of the $(n{-}1)$-dimensional permutahedron are  
indexed by the elements of $\frakS_n$ so that the $1$-skeleton
is the Hasse diagram of the weak order (see Figure~\ref{F:S4}).
Facets of the permutahedron are products of two lower dimensional
permutahedra, and the image of $\rho_\zeta$ is the set of
vertices in a facet.
Moreover, every facet arises in this way for a unique triple 
$(p,q,\zeta)$ with $p+q=n$ and $\zeta\in \Sh{(p,q)}$ 
(see~\cite[Exer.~2.9]{BLSWZ}). 
Such such a facet has {\em type} $(p,q)$. 
Figure~\ref{F:alphaperm} shows the image of
$\rho_{1324}$, a facet of the $3$-permutahedron of type $(2,2)$, and the
permutation $2431$.

The description~\eqref{E:defalpha} of $A^w_{u,v}$ (and hence of $\alpha^w_{u,v}$) 
can be interpreted as follows: 
Given $u\in\frakS_p$, $v\in\frakS_q$, and
$w\in\frakS_{p+q}$, the structure constant $\alpha^w_{u,v}$
counts the number of facets of type $(p,q)$ of the $(p{+}q{-}1)$-permutahedron 
such that the vertex $\rho_\zeta(u,v)$ is below $w$ and it is the maximum
vertex in that facet below $w$.

For instance, the facet $\rho_{1324}$ contributes to the structure constant
$\alpha^{2431}_{12,21}$ because the vertex $\rho_{1324}(12,21)=1423$ satisfies
the required properties in relation to the 
vertex $w=2431$, as shown in Figure~\ref{F:alphaperm}.

%%%%%%%%%%%%%%%%%%%%%%%%%%%%%%%%%%%%%%%%%%%%%%%%%%%%%%%%%%%%%%%%%%%%%%
\subsection{The antipode of $\SSym$}\label{S:antipode}
Malvenuto left open the problem of an explicit formula for the antipode of
$\SSym$~\cite[pp. 59-60]{Malv}. 
We identify the coefficients
of the antipode in terms of both bases in explicit combinatorial terms. 
These are based upon a general formula for the antipode of a connected Hopf
algebra due to Milnor and Moore~\cite{MM65}.

For any subset $\setS=\{p_1<p_2<\dotsb<p_k\}\subseteq[n{-}1]$ and 
$v\in\frakS_n$ set
 \[
   v_\setS\ :=\ 
      \st(v_1,\ldots,v_{p_1})\times\st(v_{p_1+1},\ldots,v_{p_2})\times\dotsb\times
       \st(v_{p_k+1},\ldots,v_{n})\ \in\ \frakS_{\setS}\,.
\]

\begin{thm}\label{T:ant-fundamental}
 For $v,w\in\frakS_n$ set
 \begin{align*}
     \lambda(v,w)\ :=\ &\ \#\{\setS\subseteq[n{-}1]\mid
     \Des(w^{-1}v_\setS)\subseteq\setS \text{ and $\#\setS$ is odd}\}\\
       & -\#\{\setS\subseteq[n{-}1]\mid \Des(w^{-1}v_\setS)\subseteq\setS 
         \text{ and $\#\setS$ is even}\}.
 \end{align*}
 Then
 \begin{equation}\label{E:ant-fundamental}
   S(\calF_v)\ =\ \sum_{w\in\frakS_n}\lambda(v,w)\,\calF_w\,.
 \end{equation}
\end{thm}

These coefficients of the antipode may indeed be
positive or negative. 
For instance,
  \[
     S(\calF_{231})\ =\ \calF_{132}-\calF_{213}-2\calF_{231}+\calF_{312}\,.
 \]
The coefficient of $\calF_{312}$ is $1$ because $\{1\}$, $\{2\}$, and
$\{1,2\}$ are the subsets $\setS$ of $\{1,2\}$ which satisfy 
$\Des\bigl((312)^{-1}(231)_\setS\bigr)\subseteq\setS$.

Our description of these coefficients is semi-combinatorial, in the sense
that it involves a difference of cardinalities of sets. 
On the monomial basis the situation is different. 
The sign of the coefficients of $S(\calM_v)$ only depends on the number of
global descents of $v$. 
We provide a fully combinatorial description of these coefficients.
Let $v,w\in\frakS_n$ and suppose $\setS\subseteq\GDes(v)$.
Define $C_\setS(v,w)\subseteq\Sh{\setS}$ to be those $\zeta\in\Sh{\setS}$
satisfying 
 \begin{equation}\label{E:C-def}
   \begin{array}{rl}
       (i) & \text{$v_\setS\zeta^{-1}\leq w$,}\\
      (ii) & \text{if $v\leq v'$ and $v'_\setS\zeta^{-1}\leq w$ then $v=v'$, \
                 and}\\
     (iii) &  \text{if $\Des(\zeta)\subseteq\setR\subseteq\setS$ and 
              $v_\setR\zeta^{-1}\leq w$ then  $\setR=\setS$.}
   \end{array}
 \end{equation}
Set $\kappa(v,w):=\#C_{\GDes(v)}(v,w)$.

\begin{thm}\label{T:ant-monomial} 
  For $v,w\in\frakS_n$, we have
 \begin{equation}\label{E:ant-monomial}
   S(\calM_v)\ =\ (-1)^{\#\GDes(v)+1}\sum_{w\in\frakS_n}\kappa(v,w)\,\calM_w\,.
 \end{equation}
\end{thm}

For instance,
 \begin{multline*}
   S(\calM_{3412})\ =\ \calM_{1234}+2\calM_{1324}+\calM_{1342}+\calM_{1423}\\
    +\calM_{2314}+\calM_{2413}+\calM_{3124}+\calM_{3142}+\calM_{3412}\,.
 \end{multline*}
Consider the coefficient of $\calM_{3412}$. 
In this case, $\setS=\GDes(3412)=\{2\}$, so
 \[
   \Sh{\setS}\ =\ \{1234,\ 1324,\ 1423,\ 2314,\ 2413,\ 3412\}\,.
 \]
We invite the reader to verify that 3412 is the only 
element of $\Sh{\{2\}}$ that satisfies all three conditions of~\eqref{E:C-def}. 
Therefore $C_{\setS}(3412,3412)=\{3412\}$ and the coefficient
is $\kappa(3412,3412)=1$.

\begin{rem}
 The antipode of $\SSym$ has infinite order. 
 A computation gives that
 \[
   S^{2m}(\calM_{231})\ =\ \calM_{231}+2m(\calM_{213}-\calM_{132})
    \qquad\forall\ m\in\Z\,.
 \]
\end{rem}

\section{Hopf-algebraic structure of $\SSym$}

\subsection{Cofreeness of $\SSym$ and the coradical filtration}
\label{S:cofree}

The basis $\{\calM_u\}$ reveals the existence of a second coalgebra grading on
$\SSym$, given by the number of global descents.
With respect to this grading, $\SSym$ is a cofree graded
coalgebra. We deduce an elegant description of the coradical
filtration: it corresponds to a filtration of the symmetric groups by certain
lower order ideals determined by the number of global descents.
In particular, the space of primitive elements is spanned by those 
$\calM_u$ where $u$ has no global descents.

Let $V$ be a vector space and consider the graded vector space
 \[ Q(V)\ :=\ \bigoplus_{k\geq 0} V^{\ten k}\,.\]
The space $Q(V)$ is a graded connected coalgebra under the
\emph{deconcatenation} coproduct
 \[
   \Delta(v_1\ten\dotsb\ten v_k)\ =\ \sum_{i=0}^k (v_1\ten\dotsb\ten v_i)\otimes
   (v_{i+1}\ten\dotsb\ten v_k)\,,
 \]
and counit $\epsilon(v_1\ten\dotsb\ten v_k)=0$  if $k\geq 1$.  

The following universal
property is satisfied for the canonical  projection $\pi:Q(V)\to V$.
Given a graded coalgebra $C=\oplus_{k\geq 0}C^k$ and a linear map $\varphi:C\to
V$ such that $\varphi(C^k)=0$ for $k\neq 1$, there is a unique morphism of
graded coalgebras 
$\hat{\varphi}:C\to Q(V)$ such that the following diagram commutes
 \[
   \xymatrix{{\ C\ }\ar@{-->}[rr]^{\hat{\varphi}}\ar[dr]_{\varphi} &
    &{Q(V)}\ar[ld]^{\pi}\\ & {V} }
\]
Explicitly, $\hat{\varphi}$  is defined by
\begin{equation} \label{E:cofree}
\hat{\varphi}_{|_{C^k}}=\varphi^{\ten k}\Delta^{(k-1)}\,.
\end{equation}
In particular, $\hat{\varphi}_{|_{C^0}}=\epsilon$,
$\hat{\varphi}_{|_{C^1}}=\varphi$ and
$\hat{\varphi}_{|_{C^2}}=(\varphi\ten\varphi)\Delta$.
We say that the graded coalgebra $Q(V)$ is cofreely cogenerated by the
projection $Q(V)\twoheadrightarrow V$.

To establish the cofreeness of $\SSym$, we first define a second coalgebra grading.
Let $\frakS^{0}:=\frakS_0$, and for
$k\geq 1$, let
 \begin{align*}
   \frakS_n^{k}&\ :=\ \{u\in\frakS_n\mid u
                \text{ has exactly $k{-}1$ global descents}\}, \
       \text{ and }\\
    \frakS^{k}&\ :=\ \coprod_{n\geq 0}\frakS_n^{k}\,.
 \end{align*}
For instance,
 \begin{multline*}
  \frakS^{1}\ =\ \{1\}\ \cup\ \{12\}\ \cup\ \{123,\,213,\,132\}\
                 \cup\ \{1234,\,2134,\,1324,\,1243,\,3124,\\
                   2314,\,2143,\,1423,\,1342,\,3214,\,
                   3142,\,2413,\,1432\}\ \cup\ \dotsb
 \end{multline*}

Let $(\SSym)^k$ be the vector subspace of $\SSym$ spanned by
 $\{\calM_u\mid u\in\frakS^k\}$.

\begin{thm}\label{T:cofree} The decomposition $\SSym=\oplus_{k\geq 0}(\SSym)^k$
is a coalgebra grading. Moreover, endowed with this grading, $\SSym$ is a
cofree graded coalgebra. \end{thm}

The coradical $C^{(0)}$ of a graded connected coalgebra $C$ is the 1-dimensional
component in degree 0 (identified with the base field via the counit).
The primitive elements %of $C$
 are
 \[
   \text{Prim}(C)\ :=\ \{x\in C\mid \Delta(x)=x\ten 1+1\ten x\}\,.
 \]
Set $C^{(1)}:=C^{(0)}\oplus \text{Prim}(C)$, the first level of the
coradical filtration.
More generally, the $k$-th level of the coradical filtration is
 \[
   C^{(k)}\ :=\ \bigl(\Delta^{(k)}\bigr)^{-1}
        \Bigl(\sum_{i+j=k}C^{\ten i}\ten C^{(0)}\ten C^{\ten j}\Bigl)\,.
 \]
We have $C^{(0)}\subseteq C^{(1)}\subseteq C^{(2)}
          \subseteq\dotsb\subseteq C=\bigcup_{k\geq 0}C^{(k)}$,
and
 \[
   \Delta(C^{(k)})\ \subseteq\ \sum_{i+j=k}C^{(i)}\ten C^{(j)}\,.
 \]
Thus, the coradical filtration measures the complexity of
iterated coproducts.

For a cofree graded coalgebra $Q(V)$, the coradical filtration is easy to
describe. The space of primitive elements is just $V$, and the $k$-th level
of the coradical filtration is $\oplus_{i=0}^k V^{\ten i}$. 
These are immediate from the definition of the deconcatenation
coproduct.

Define
\[ \frakS_n^{(k)}\ :=\ \coprod_{i=0}^k \frakS_n^k \text{ \ and \ }
 \frakS^{(k)}\ :=\ \coprod_{i=0}^k \frakS^k \,.\]
 In other words,  $\frakS^{(0)}=\frakS_0$ and for
$k\geq 1$,
\[\frakS_n^{(k)}=\{u\in\frakS_n\mid u
                \text{ has at most $k{-}1$ global descents}\}\,.\]

Proposition~\ref{P:galoisglobal} asserts that
$\GDes\colon\frakS_n\to\calQ_n$ is order-preserving.
Since $\calQ_n$ is ranked by the cardinality of a subset, it follows that
  $\frakS_n^{(k)}$ is a lower order ideal of $\frakS_n$, with
  $\frakS_n^{(k)}\subseteq\frakS_n^{(k+1)}$.
The coradical filtration corresponds precisely to this filtration of the
symmetric groups by lower ideals.

 \begin{coro}\label{C:coradical}
 A linear basis for the $k$-th level of the coradical filtration of $\SSym$ is
 \[
   \{\calM_u\mid u\in \frakS^{(k)}\}\,.
 \]
 In particular, a linear basis for the space of primitive elements is
 \[
    \{\calM_u\mid u \text{ has no global descents}\}\,.
 \]
\end{coro}

%%%%%%%%%%%%%%%%%%%%%%%%%%%%%%%%%%%%%%%%%%%%%%%%%%%%%%%%%%%%%%%%%%%%%%%%%
\subsection{The descent map to quasi-symmetric functions}\label{S:descentmap}

We study the effect of the morphism of Hopf algebras~\eqref{E:descentmap}
 \[
   \calD\ \colon\ \SSym\ \to\ \QSym, \quad\text{ defined by }\quad
   \calF_u\ \mapsto\  F_{\Des(u)}\,
 \]
on the monomial basis.
Here, we use subsets $\setS$ of $[n{-}1]$ to index monomial 
quasi-symmetric functions of degree $n$.

\begin{defi}\label{D:closed}
A permutation $u\in\frakS_n$ is {\em closed} if we have 
$u=Z(\setT)$ for some $\setT\in\calQ_n$.
Equivalently, $u$ is
closed if and only if $\Des(u)=\GDes(u)$.
\end{defi}

\begin{thm}\label{T:map-monomial}
 Let $u\in\frakS_n$. Then
 \[
   \calD(\calM_u)\ =\ \begin{cases}
           M_{\GDes(u)} & \text{if $u$ is closed,}\\
                     0 & \text{if not.}
        \end{cases}
 \]
\end{thm}

\subsection{$\SSym$ as a crossed product over $\QSym$}\label{S:crossed}
We describe the {\it algebra} 
structure of $\SSym$ as a crossed product over the Hopf algebra $\QSym$.
See~\cite[\S 7]{Mo93a} for a review of this construction
in the general Hopf algebraic setting.  Let us only say that the crossed
product of a Hopf algebra $K$ with an algebra $A$ with respect to a Hopf cocycle
$\sigma:K\ten K\to A$ is a certain algebra structure on the space $A\otimes K$,
denoted by $A\#_{\sigma} K$.

\begin{prop}
 The map $M_{\setS}\mapsto \calM_{Z(\setS)}$
 induces a morphism of  coalgebras $\calZ\colon\QSym\to\SSym$ that is a right inverse
 to the morphism of Hopf algebras $\calD\colon\SSym\to\QSym$.
\end{prop}

In this situation, an important theorem of Blattner, Cohen and
Montgomery~\cite{BCM86} applies. Namely, suppose $\pi:H\to K$ is a morphism of
Hopf algebras that admits a coalgebra splitting $\gamma:K\to H$. Then
there is a \emph{crossed product} decomposition
\[H\ \cong\ A\#_{\sigma}K\]
where $A$ is the \emph{left Hopf kernel} of $\pi$:
\[A\ =\ \{h\in H \mid \sum h_1\ten\pi(h_2)=h\ten 1\}\]
and the \emph{Hopf cocycle} $\sigma:K\ten K\to A$ is
\[\sigma(k,k')=\sum\gamma(k_1)\gamma(k'_1)S\gamma(k_2k'_2)\,.\]
Note that if $\pi$ and $\gamma$ preserve
gradings, then so does the rest of the structure.

Let $A$ be the left Hopf kernel of $\calD\colon\SSym\to\QSym$ and
$A_n$ its $n$-th homogeneous component.
Once again the monomial basis of $\SSym$ proves useful in
describing $A$.
\begin{prop} A basis for  $A_n$ is the
set  $\{\calM_u\}$ where $u$ runs over all permutations of $n$ that
are not of the form
 \[
    *\ldots*12\ldots n{-}k
 \]
for any $k=0,\ldots,n-1$. In particular,
\[\dim A_n\ =\ n!-\sum_{k=0}^{n-1}k!\,.\]
\end{prop}
%%%%%%%%%%%%%%%%%%%%%%%%%%%%%%%%%%%%%%%%%%%%%%%%%%%%%%%%%%%%%%%%%%%%%%

\def\cprime{$'$}
\providecommand{\bysame}{\leavevmode\hbox to3em{\hrulefill}\thinspace}
\providecommand{\MR}{\relax\ifhmode\unskip\space\fi MR }
\providecommand{\MRhref}[2]{%
  \href{http://www.ams.org/mathscinet-getitem?mr=#1}{#2}
}
\providecommand{\href}[2]{#2}

\end{document}